\documentclass[11pt,a4paper]{article}

\usepackage{latexsym}
\usepackage{amsmath, amsthm, amsfonts, amssymb}
\usepackage{epsfig} 
\usepackage{hyperref}

\newtheorem{thr}{Theorem}
\newtheorem{lem}{Lemma}
\newtheorem{cor}{Corollary}

\begin{document}

%\pagenumbering{gobble}
%\maketitle
\centerline{\Large{\bf Some properties of generalized self-reciprocal}}
\centerline{}
\centerline{\Large{\bf polynomials over finite fields}}
\centerline{}
\centerline{\textsuperscript{a,*}Ryul Kim, \textsuperscript{b}Ok-Hyon Song, \textsuperscript{c}Hyon-Chol Ri}
\centerline{}
\footnotesize \centerline{\textsuperscript{a,b}Faculty of Mathematics, \textbf{Kim Il Sung} University, DPR Korea}
\footnotesize \centerline{\textsuperscript{c}Department of Mathematics, Chongjin University of Education No.2, DPR Korea}
\footnotesize \centerline{\textsuperscript{*}Corresponding author. e-mail address: ryul\_kim@yahoo.com}
\centerline{}
\centerline{}
%\hline
\begin{abstract}
Numerous results on self-reciprocal polynomials over finite fields have been studied. In this paper we generalize some of these to $a$-self reciprocal polynomials defined in \cite{fit}. We consider some properties of the divisibility of $a$-reciprocal polynomials and characterize 
the parity of the number of irreducible factors for $a$-self reciprocal polynomials over finite fields of odd characteristic.

\end{abstract}


%\hline
%
%-----------------------1. Introduction-----------------
%                                                      
\section{Introduction}

Self-reciprocal polynomials over finite fields are of interest both from a theoretical and a practical viewpoint, 
so they are widely studied by many authors (see, for example, \cite{ahm2, car, gar1}).
Carlitz \cite{car} proposed a formula on the number of self-reciprocal irreducible monic (srim) polynomials over finite fields and 
Meyn \cite{mey} gave a simpler proof of it. Recently Carlitz's result was generalized by Ahmadi \cite{ahm1} and similar explicit results have been obtained in \cite{gar2, stich}. Using the Stickelberger-Swan theorem, Ahmadi and Vega \cite{ahm2} characterized the parity of the number of irreducible factors of a self-reciprocal polynomials over finite fields. 
In \cite{yuc}, Yucas and Mullen classified srim polynomials based on their orders and considered the weight of srim polynomials. 
The problem concerning to the existence of srim polynomials with prescribed coefficients has also been considered, see \cite{gar1}. 

On the other hand, Fitzgerald and Yucas \cite{fit} introduced a new notion of generalized reciprocal polynomial over finite fields to give new descriptions of the factors of Dickson polynomials over finite fields. They characterized the generalized reciprocal polynomials by their roots and orders so that the results of \cite{yuc} were generalized. 
The generalization of reciprocal polynomials gives a possibility to find new special types of irreducible polynomials over finite fields. 
This motivated us to consider the generalized reciprocal polynomials over finite fields.
Consideration of the divisibility and the number of irreducible factors of a given polynomial is the first step to test its irreducibility.

In this work we find some properties for the divisibility of generalized reciprocal and generalized self-reciprocal polynomials for the 
purpose of characterizing a class of the generalized srim polynomials. 
First we characterize the generalized reciprocal of the product of given polynomials, extend the notion of generalized self-reciprocal 
polynomial defined only for even degree in \cite{fit} to odd degree and obtain some new results on the divisibility of 
trivial or nontrivial generalized self-reciprocal polynomials. 
Next we determine the parity of the number of irreducible factors of generalized self-reciprocal polynomials which is a generalization of the results in \cite{ahm2}.
In fact, we also have obtained an explicit formula on the number generalized srim polynomials but when we prepare for submission, 
a more generalized result in \cite{ahm1} came to our view, so the formula was omitted from this paper.  

Throughout the paper, let $q=p^e$ be an odd prime power, $\mathbf{F}_q$ be a finite field containing $q$ elements and $0 \neq a \in \mathbf{F}_q$.  
For $f(x) \in \mathbf{F}_q[x]$, a monic polynomial of degree $n$ with $f(0) \neq 0$, $\hat{f}_a(x) := \frac{x^n}{f(0)}f \left( \frac{a}{x} \right)$ 
is called \textsl{$a$-reciprocal} of $f(x)$ \cite{fit}. 
That is, if $f(x)=\sum^n_{\substack{i=0}} b_i x^i$ then $\hat{f}_a(x)=(1/b_0)\sum^n_{\substack{i=0}} b_i a^ix^{n-i}$. $\hat{f}_a(x)$ is also monic and if $\alpha$ is a root of $f(x)$ then $a/\alpha$ is a root of $\hat{f}_a(x)$. 
The case $a=1$ is the usual reciprocal in \cite{ahm2,mey,yuc}.

%
%------2. Divisibility of generalized reciprocal polynomials -----------------
%
\section{Divisibility of generalized reciprocal polynomials}

In this section we summarize some properties for the divisibility of generalized reciprocal and generalized self-reciprocal polynomials over $\mathbf{F}_q$.

%
%-------------Theorem 1---------------
%

\begin{thr}\label{theo1}
$a$-reciprocal of a product of two polynomials is the product of $a$-reciprocals of given polynomials.
\end{thr}
\begin{proof}
Let $f(x)=\sum^n_{\substack{i=0}} b_i x^i, ~ g(x)=\sum^m_{\substack{j=0}} c_j x^j , (b_0 \neq 0, c_0 \neq 0)$. From the definition,
\begin{eqnarray*}
\hat{f}_a(x)=\frac{1}{b_0} \sum^n_{\substack{i=0}} b_i a^i x^{n-i} = \frac{1}{b_0} \sum^n_{\substack{i=0}} b_{n-i} a^{n-i} x^i, \\
\hat{g}_a(x)=\frac{1}{c_0} \sum^m_{\substack{j=0}} c_j a^j x^{m-j} = \frac{1}{c_0} \sum^m_{\substack{j=0}}c_{m-j} a^{m-j} x^j,
\end{eqnarray*}
\begin{align*}
\hat{f}_a(x) \cdot \hat{g}_a(x) &=\frac{1}{b_0 c_0} \sum^{n+m}_{\substack{k=0}} \left[ \sum_{\substack{i+j=k}} \left( b_{n-i} a^{n-i} \right) \left( c_{m-j} a^{m-j} \right) \right] x^k \\
&= \frac{1}{b_0 c_0} \sum^{n+m}_{\substack{k=0}} \left( \sum_{\substack{i+j=k}} b_{n-i} \cdot c_{m-j} \cdot a^{n+m-k} \right) x^k
\end{align*}
and
\begin{equation*}
f(x) \cdot g(x) = \sum^{n+m}_{\substack{k=0}} \left( \sum_{\substack{i+j=k}} b_i c_j \right) x^k = \sum^{n+m}_{\substack{k=0}} d_k x^k \left( d_k := \sum_{\substack{i+j=k}} b_i c_j \right) .
\end{equation*}
Therefore
\begin{equation*}
\widehat{fg}_a(x)=\frac{1}{b_0 c_0} \sum^{n+m}_{\substack{k=0}} d_k a^k x^{n+m-k} = \frac{1}{b_0 c_0} \sum^{n+m}_{\substack{k=0}} d_{n+m-k} a^{n+m-k} x^k
\end{equation*}
We set $i':=n-i, j':=m-j$, then
\begin{equation*}
d_{n+m-k}= \sum_{\substack{i'+j'=k}} b_{n-i'} c_{m-j'}
\end{equation*}
and thus
\begin{equation*}
\widehat{fg}_a(x) = \frac{1}{b_0 c_0} \sum^{n+m}_{\substack{k=0}} \left( \sum_{\substack{i'+j'=k}} b_{n-i'} \cdot c_{m-j'} \cdot a^{n+m-k} \right) x^k.
\end{equation*}
\end{proof}

It might be $\hat{f}_a (x)=\hat{f}_b (x)$ for nonzero distinct elements $a, b$ in $\mathbf{F}_q$. For example, if $f(x)=x^2+c \in \mathbf{F}_5[x]$, $c \neq 0$, then $\hat{f}_2 (x)=\hat{f}_3 (x)$.

A monic polynomial $f(x) \in \mathbf{F}_q [x]$ is said to be \textsl{$a$-self reciprocal} when $\hat{f}_a (x)=f (x)$.
In \cite {fit}, the notion of $a$-self reciprocal polynomial has been defined for only even degree. 
Note that $f(x)= \sum^n_{\substack{i=0}} b_i x^i$ is $a$-self reciprocal if and only if for each $i$, $b_{n-i} b_0=b_i a^i$. When $i=n$, we see that $b^2_0=a^n$.

First consider a case when $n$ is odd. If $a^n$ and thus $a$ is not a square in $\mathbf{F}_q$, then there is not any $a$-self reciprocal monic (srm) polynomial of degree $n$. When $0 \neq a \in \mathbf{F}_q$ is a square in $\mathbf{F}_q$, if $f(x)$ is $a$-srm polynomial then either $b_0 = \sqrt{a}^n$ or $b_0 = -\sqrt{a}^n$. 
%
%-------------Theorem 2---------------
%
\begin{thr}\label{theo2}
Let $n$ be an odd, $0 \neq a \in \mathbf{F}_q$ be a square in $\mathbf{F}_q$ and $f(x) \in \mathbf{F}_q[x]$ be $a$-srm polynomial of degree $n$.

$(1)$ If $b_0=\sqrt{a}^n$ then $f(x)$ is divided by $x+\sqrt{a}$.

$(2)$ If $b_0=-\sqrt{a}^n$ then $f(x)$ is divided by $x-\sqrt{a}$.
\end{thr}
\begin{proof}
Suppose $b_0=\sqrt{a}^n$. Then we have
\begin{align*}
f \left( -\sqrt{a} \right) &=\left( -\sqrt{a} \right) ^n+b_{n-1} \left( -\sqrt{a} \right) ^{n-1}+ \cdots +b_1\left( -\sqrt{a} \right) + \left( \sqrt{a} \right) ^n \\
&= b_{n-1} \left( -\sqrt{a} \right) ^{n-1}+ \cdots +b_1\left( -\sqrt{a} \right) ,\\
\hat{f}_a\left( -\sqrt{a} \right) &=- \frac{1}{\left( -\sqrt{a} \right) ^n} \left[ \sqrt{a}^n \left( -\sqrt{a} \right) ^n+%b_1 a \left( -\sqrt{a} \right) ^{n-1}+ %
\cdots +b_{n-1} a^{n-1} \left( -\sqrt{a} \right) +a^n \right] \\
&= -\left[ \sqrt{a}^n+b_1 \left( -\sqrt{a} \right) + \cdots +b_{n-1} \left( -\sqrt{a} \right)^{n-1}+\left( -\sqrt{a} \right)^n \right] \\
&= -\left[ b_1 \left( -\sqrt{a} \right) + \cdots +b_{n-1} \left( -\sqrt{a} \right)^{n-1} \right] = -f \left( -\sqrt{a} \right)
\end{align*}
and from $f=\hat{f}_a$, $f \left( -\sqrt{a} \right)= -f \left( -\sqrt{a} \right)$ which implies that $f \left( -\sqrt{a} \right)= 0$. The second case follows similarly. 
\end{proof}

%
%-------------Corrollary 1---------------
%
\begin{cor}\label{cor1}
If $0 \neq a \in \mathbf{F}_q$ is a square in $\mathbf{F}_q$, then $x+\sqrt{a}$ and $x-\sqrt{a}$ are the only $a$-srim polynomials of odd degree over $\mathbf{F}_q$.
\end{cor}

Next let $n$ be an even, that is, $n=2m$. $f(x)= \sum^{2m}_{\substack{i=0}} b_i x^i$ is said to be \textsl{trivial} or \textsl{nontrivial} respectively, according to $b_0=-a^m$ or $b_0=a^m$\cite {fit}.
If $f(x)$ is a trivial $a$-srm polynomial, then $f \left( \sqrt{a} \right) = f \left( -\sqrt{a} \right) =0$, hence $f(x)$ is a multiple of $x^2-a$. Therefore $x^2-a$ is the only trivial $a$-srim polynomial. From the definition, $b_m b_0 = b_m a^m$ and $b_0 = -a^m$, so that $b_m = 0$, hence trivial $a$-srm polynomials have always even terms.

%
%-------------Lemma 1---------------
%
\begin{lem}\label{lem1}
$x^2-a$ is a trivial $a$-srm polynomial and $\left( x^2-a \right) ^2$ is a nontrivial one.
\end{lem}

%
%-------------Lemma 2---------------
%

\begin{lem}\label{lem2}
The product of two $a$-srm polynomials satisfies the following multiplication table.

\bigskip

\begin{tabular}[!hbp]{c|c c}
$\cdot$ & trivial & nontrivial \\
\hline
trivial & nontrivial & trivial \\  

nontrivial & trivial & nontrivial \\  

\end{tabular}

\bigskip

\end{lem}
The proof of the above lemmas are very simple so we omit.

%
%--------------Theorem 3----------------
%
\begin{thr}\label{theo3}
If $f(x) \in \mathbf{F}_q[x]$ is $a$-srm polynomial of even degree, then it can be written as $f(x)=\left( x^2-a \right) ^k \cdot g(x)$, 
where $g(x)$ is a nontrivial $a$-srm polynomial not divided by $x^2-a$ and

\begin{equation*}
k = \left\{ 
\begin {array}{ll}
\textrm{odd}, & f(x) \textrm{ is trivial},\\
\textrm{even}, & f(x) \textrm{ is nontrivial} 
\end{array} \right..
\end{equation*}
\end{thr}
\begin{proof}
Suppose that $x^2-a$ divides $f(x)$. Then $f(x)=\left( x^2-a \right) \cdot f_1(x)$ for some $f_1(x) \in \mathbf{F}_q[x]$. 
If $f(x)$ is trivial, then $f_1(x)$ is nontrivial by Theorem 1 and the above lemmas. And if $f(x)$ is nontrivial, then $f_1(x)$ is trivial, 
so it is divided by $x^2-a$ and it can be written as $f(x)=\left( x^2-a \right) ^2 \cdot f_2(x)$. Again from the above lemmas $f_2(x)$ is nontrivial. Continuing this procedure implies the claim. 
\end{proof}
 
If $f(x) \in \mathbf{F}_q[x]$ is divided by $x^2-a$, then $\sqrt{a}$ and $-\sqrt{a}$ are roots of $f(x)$. If $a$ is a square in $\mathbf{F}_q[x]$, then there might be an $a$-srm polynomial which is not divided by $x^2-a$ but has a root $\sqrt{a}$ or $-\sqrt{a}$.
For example, if $q=5, a=4$, then $\sqrt{a}=2, -\sqrt{a}=3$,
\begin{equation*}
f(x)=x^6+x^5+3x^4+4x^3+2x^2+x+4 \in \mathbf{F}_5[x]
\end{equation*}
is a nontrivial $4$-srm polynomial which is not divided by $x^2-4$ but $f(2)=0$ and
\begin{equation*}
g(x)=x^6+3x^4+3x^3+2x^2+4 \in \mathbf{F}_5[x]
\end{equation*}
is a nontrivial $4$-srm polynomial which is not divided by $x^2-4$ but $g(3)=0$.

%
%---------------------Theorem 4-----------------------
%
\begin{thr}\label{theo4}
Let $f(x) \in \mathbf{F}_q[x]$ be a nontrivial $a$-srm polynomial of degree $n=2m$ and $a$ be a square in $\mathbf{F}_q$. If $f(x)$ is not divided by $x^2-a$ and $f \left( \sqrt{a} \right)=0$, then it can be written as $f(x)=\left( x- \sqrt{a} \right) ^k \cdot g(x)$, where $k$ is even and $g(x)$ is a nontrivial $a$-srm polynomial with $g \left( \sqrt{a} \right) \neq 0$.
\end{thr}	
\begin{proof}
Let $f(x)= \sum^n_{\substack{i=0}} b_i x^i$ then by assumption $b_0=f(0)=a^m = \sqrt{a} ^n$. Set
\begin{equation*}
f_1(x):=\frac {f(x)}{x- \sqrt{a}} = x^{n-1}+c_{n-2}x^{n-2}+ \cdots +c_1x+c_0 \in \mathbf{F}_q[x].
\end{equation*}
Clearly $\left( -\sqrt{a} \right) \cdot c_0 = b_0 = \sqrt{a} ^n$ which implies that $c_0 = -\sqrt{a} ^{n-1}$. Then $f_1(x)$ is $a$-srm polynomial since $x-\sqrt{a}$ is $a$-srm polynomial. By Theorem 2, $f_1(x)$ is divided by $x-\sqrt{a}$ and thus we can write as $f(x)=\left( x- \sqrt{a} \right) ^2\cdot f_2(x)$. It is clear that $\left( x- \sqrt{a} \right) ^2$ is a nontrivial $a$-srm polynomial, hence $f_2(x)$ is a nontrivial $a$-srm polynomial of degree $n-2$. Continuing this procedure for $f_2(x)$ completes the proof. 
\end{proof}

When $a$ is a square in $\mathbf{F}_q$, from this theorem the consideration of nontrivial $a$-srm polynomials not divided by $x^2-a$ is reduced to the case of ones without roots $\sqrt{a}$ and $-\sqrt{a}$.

%
%------3. The parity of the number of irreducible factors of generalized self-reciprocal polynomials-----------------
%
\section{The parity of the number of irreducible factors of generalized self-reciprocal polynomials}

In this section we characterize the parity of the number of irreducible factors of the nontrivial $a$-srm polynomials over $\mathbf{F}_q$.

%---------------------Theorem 5-----------------------
%
\begin{thr}\label{theo5}
Let $f(x) \in \mathbf{F}_q[x]$ be a nontrivial $a$-srm polynomial of degree $2n$ with $r$ pairwise distinct irreducible factors over $\mathbf{F}_q$. Then $r \equiv 0 \pmod{2}$ if and only if $(-1)^n a^{n(n-2)} f \left( \sqrt{a} \right) f \left( -\sqrt{a} \right) $ is a square in $\mathbf{F}_q$, where $\sqrt{a}$ is a square root of $a$ in an extension of $\mathbf{F}_q$. 
\end{thr}
\begin{proof}
Let $f(x)= \sum^{2n}_{\substack{i=0}} b_i x^i$, then $f(x)= x^n g \left( x+ \frac{a}{x} \right)$, where $g(x)= b_n + \sum^{n-1}_{\substack{i=0}}$ $ b_{2n-i}D_{n-i, a}(x)$ and $D_{n-i, a}(x)$ is the Dickson polynomial $\cite{fit}$.
Clearly $f(0)= a^n \neq 0$ and by properties of the resultant, $R(f, f') = R(f, nf-xf')$. Since
\begin{align*}
f'(x) &= nx^{n-1} g(x+a/x)+x^n g'\left(x+\frac{a}{x} \right) \left(1-\frac{a}{x^2} \right)\\
&= nx^{n-1} g\left(x+\frac{a}{x} \right)+x^{n-2} g'\left(x+\frac{a}{x} \right)(x^2-a)
\end{align*}
and 
\begin{equation*}
nf(x)-xf'(x) = -x^{n-1} g'\left(x+\frac{a}{x} \right)(x^2-a),
\end{equation*}
we have
\begin{align*}
D(f) &= (-1)^{n(2n-1)} R(f, f') = (-1)^n \cdot f(0)^{-1} \cdot R(f, nf-xf') \\
&= (-1)^n \cdot a^{-n} \cdot R(f, x^2-a) \cdot R\left(f, -x^{n-1} g'\left(x+\frac{a}{x} \right)\right) \\
&= (-1)^n \cdot a^{-n} \cdot f \left( \sqrt{a} \right) f \left( -\sqrt{a} \right) \cdot R\left(f, -x^{n-1} g'\left(x+\frac{a}{x} \right)\right)  
\end{align*}
Let $x_0, \cdots, x_{n-1}, \frac{a}{x_0}, \cdots, \frac{a}{x_{n-1}}$ be the roots of $f(x)$ in some extension of $\mathbf{F}_q$. Then
\begin{align*}
R \left( f, -x^{n-1} g' \left( x+ \frac{a}{x} \right) \right) &= \prod^{n-1}_{\substack{i=0}} x_i^{n-1} g' \left( x_i+\frac{a}{x_i} \right) \prod^{n-1}_{\substack{i=0}} \left( \frac{a}{x_i} \right) ^{n-1} g' \left( \frac{a}{x_i}+x_i \right) \\
&= \prod^{n-1}_{\substack{i=0}} x_i^{n-1} g' \left( x_i+\frac{a}{x_i} \right) \prod^{n-1}_{\substack{i=0}} a^{n-1} x_i^{1-n} g' \left( x_i+\frac{a}{x_i} \right) \\
&= a^{n(n-1)} \left[ \prod^{n-1}_{\substack{i=0}} g' \left( x_i+\frac{a}{x_i} \right) \right] ^2.
\end{align*}
On the other hand, since $x_i$ is a root of $f(x)$ if and only if $x_i + \frac{a}{x_i}$ is a root of $g(x)$ and
\begin{equation*}
D(g) = (-1)^{n(n-1)/2} R(g, g') = (-1)^{n(n-1)/2} \prod^{n-1}_{\substack{i=0}} g' \left( x_i+\frac{a}{x_i} \right),
\end{equation*}
we have
\begin{equation*}
D(f) = (-1)^n a^{n(n-2)} f \left( \sqrt {a} \right) f \left( -\sqrt {a} \right) \cdot D(g)^2.
\end{equation*}
Required result follows from  Stickelberger theorem $\cite{stick}$ (see also [2, Theorem 4]). 
\end{proof}

Let $f(x)= \sum^{2n}_{\substack{i=0}} b_i x^i$ be as in Theorem 5, then
\begin{equation*}
f(x)= \sum^n_{\substack{i=0}} b_{2i} x^{2i} + \sum^{n-1}_{\substack{i=0}} b_{2i+1} x^{2i+1}
\end{equation*}
and
\begin{eqnarray*}
f \left( \sqrt{a} \right) = \sum^n_{\substack{i=0}} b_{2i} a^i + \sqrt{a} \sum^{n-1}_{\substack{i=0}} b_{2i+1} a^i, \\
f \left( -\sqrt{a} \right) = \sum^n_{\substack{i=0}} b_{2i} a^i - \sqrt{a} \sum^{n-1}_{\substack{i=0}} b_{2i+1} a^i.
\end{eqnarray*}
Denote
\begin{equation*}
A := \sum^n_{\substack{i=0}} b_{2i} a^i, \qquad B := \sum^{n-1}_{\substack{i=0}} b_{2i+1} a^i
\end{equation*}
then $f \left( \sqrt{a} \right) f \left( -\sqrt{a} \right) = A^2-aB^2$, which means $f \left( \sqrt{a} \right) f \left( -\sqrt{a} \right) \in \mathbf{F}_q$.
Since $f(x)$ is nontrivial, $b_i = b_{2n-i} a^{n-i}$. Therefore if $n$ is odd, then
\begin{equation*}
A =2\sum^{(n-1)/2}_{\substack{i=0}} b_{2n-2i} a^{n-i}, \quad B = 2\sum^{(n-3)/2}_{\substack{i=0}} b_{2n-2i-1} a^{n-i-1} +b_n a^{(n-1)/2},
\end{equation*}
and if $n$ is even then
\begin{equation*}
A =2\sum^{n/2-1}_{\substack{i=0}} b_{2n-2i} a^{n-i}+b_n a^{n/2}, \quad B = 2\sum^{n/2-1}_{\substack{i=0}} b_{2n-2i-1} a^{n-i-1}.
\end{equation*}

For any monic polynomial $f(x) \in \mathbf{F}_q[x]$ of degree $n$, denote $f^Q_a(x):=x^n f\left(x+\frac{a}{x}\right)$. Then $f^Q_a(x)$ is a nontrivial $a$-srm polynomial of degree $2n$. This is a natural generalization of $f^Q(x)$ in $\cite{mey}$ and it has been mentioned in $\cite{fit}$, too.

%---------------------Theorem 6-----------------------
%
\begin{thr}\label{theo6}
Suppose that $f(x) \in \mathbf{F}_q[x]$ is a monic irreducible polynomial of degree $n(\geq 1)$ and $f \left( \sqrt{a} \right) f \left( -\sqrt{a} \right) \neq 0$. Then $f^Q_a(x)$ is either $a$-srim polynomial or a product of $a$-reciprocal pair of irreducible polynomials of degree $n$ which are not $a$-self reciprocal.
\end{thr}
\begin{proof}
Let $\alpha$ be any root of $f^Q_a(x)$. Then $\beta = \alpha +\frac{a}{\alpha}$ is a root of $f(x)$, so $\beta ^{q^n} = \beta$. Here $n$ is the smallest positive integer $k$ with $\beta ^{q^k} = \beta$. Multiply $\alpha ^{q^n}$ to both sides of the identity
\begin{equation*}
\alpha ^{q^n} + \frac{a}{\alpha ^{q^n}} = \beta ^{q^n} = \beta = \alpha + \frac{a}{\alpha},
\end{equation*}
then
\begin{equation*}
\alpha ^{2q^n} + a = \alpha ^{q^n+1} + a \alpha ^{q^n-1}
\end{equation*}
and so
\begin{equation*}
\left( \alpha ^{q^n+1} - a \right) \left( \alpha ^{q^n-1} - 1 \right) = 0.
\end{equation*}
Hence $\alpha ^{q^n+1} = a$ or $\alpha ^{q^n-1} = 1$. When $\alpha ^{q^n+1} = a$, the irreducible factor $g(x)$ of $f^Q_a(x)$ with a root $\alpha$ has degree $\geq 2$ because if $g(x)$ has degree 1, then $a$ is a square in $\mathbf{F}_q$ and $g(x)=x \pm \sqrt{a}$  which contradicts to $f \left( \sqrt{a} \right) f \left( -\sqrt{a} \right) \neq 0$. Therefore $g(x)$ is a nontrivial $a$-srm polynomial of degree $2d$ for a positive divisor of $n$. If we assume $d<n$ then $\alpha ^{q^d+1} = a$, so $\beta ^{q^d} = \beta$, a contradiction. Hence $d=n$ and $f^Q_a(x)$ is irreducible over $\mathbf{F}_q$. When $\alpha ^{q^n-1} = 1$, the irreducible factor $g(x)$ of $f^Q_a(x)$ with a root $\alpha$ divides $x ^{q^n} - x$ and the degree of $g(x)$ is $n$. Therefore $f^Q_a(x)$ is a product of two $n$-degree irreducible polynomials $g(x)$ and $h(x) = \frac{f^Q_a(x)}{g(x)}$. Suppose $g(x)$ is a nontrivial $a$-srim polynomial, then $\alpha ^{q^{n/2}+1} - a = 0$ which contradicts to the minimality of $n$. So $g(x)$ and $h(x)$ are not $a$-self reciprocal. And if $f^Q_a(x)$ has a root $\alpha$ then it also has a root $\frac{a}{\alpha}$, which implies easily that $h(x)$ is $a$-reciprocal of $g(x)$. 
\end{proof}

Finally we prove that Theorem 5 is still valid when $f(x)$ has repeated irreducible factors.

%---------------------Theorem 7-----------------------
%
\begin{thr}\label{the7}
Let $f(x) \in \mathbf{F}_q[x]$ be a nontrivial $a$-srm polynomial of degree $2n$ with $f \left( \sqrt{a} \right) f \left( -\sqrt{a} \right) \neq 0$, and let $r$ be the number of irreducible factors counted with multiplicity of $f(x)$. Then $r$ is even if and only if $(-1)^n a^{n(n-2)}f \left( \sqrt{a} \right) f \left( -\sqrt{a} \right)$ is a square in $\mathbf{F}_q$.
\end{thr}
\begin{proof}
As in the proof of Theorem 5, there exists $g(x) \in \mathbf{F}_q[x]$ such that $f(x) = x^n g\left(x+\frac{a}{x}\right)$. Suppose that $g(x) = g_1(x) \cdots g_k(x)$ where $g_i(x) \in \mathbf{F}_q[x]$ is monic irreducible polynomial of degree $n_i$ and $n_1+ \cdots +n_k = n$. Denote $f_i(x) := x^{n_i} g_i\left(x+\frac{a}{x}\right)$, then $f(x) = f_1(x) \cdots f_k(x)$ where every $f_i(x)$ is a nontrivial $a$-srm polynomial of degree $2n_i$ over $\mathbf{F}_q$. By Theorem 6, every $f_i(x)$ is either irreducible or a product of two distinct monic irreducible polynomials of degree $n_i$. Now suppose that the claim holds for two nontrivial $a$-srm polynomials $g(x)$, $h(x)$ of degree $2s$, $2t$ with $r_s$, $r_h$ irreducible factors counted with multiplicity respectively over $\mathbf{F}_q$. If both $r_s$ and $r_h$ are odd, then neither $(-1)^s a^{s(s-2)}g(\sqrt{a}) g( -\sqrt{a})$ nor $(-1)^t a^{t(t-2)}h (\sqrt{a}) h( -\sqrt{a})$ is a square in $\mathbf{F}_q$, hence $(-1)^{s+t} a^{(s+t)(s+t-2)}g (\sqrt{a})$ $g(-\sqrt{a}) h( \sqrt{a}) h( -\sqrt{a})$ is a square in $\mathbf{F}_q$.
Similar observation about other cases for $r_s$ and $r_h$ tells us the claim is still valid for $g(x)h(x)$. It is sufficient to apply Theorem 5 to complete the proof. 
\end{proof}

{\bf Acknowledgement}. We would like to thank anonymous referees for their valuable comments and suggestions.

\end{document}